\documentclass[a4paper,oneside,portrait,12pt]{AmsArt}

\usepackage[margin=1in]{geometry}

\usepackage {graphicx}
\usepackage{amsfonts}
\usepackage{amsthm}
\usepackage{amsmath}
\usepackage{amsfonts}
\usepackage{latexsym}
\usepackage{amssymb}
\usepackage{epsfig,color}
\usepackage{graphicx, amssymb}

\usepackage{hyperref}
\hypersetup{
  bookmarks=true,
}

\newtheorem{conjecture}[]{Conjecture}

\newtheorem{theorem}{Theorem}[]

\theoremstyle{definition}

\title{Higgs bundles and applications}
\author{Laura P. Schaposnik }

\address{Department of Mathematics, University of Illinois, Chicago, IL 60607, USA.}
\email{schapos@uic.edu}

\date{\today}

\begin{document}

\markboth{Laura P. Schaposnik}
{Higgs bundles and applications}





\maketitle

\begin{abstract}
 This short  note gives an overview of how a few conjectures and theorems of the author and collaborators fit together. It was prepared for   Oberwolfach's workshop {\it Differentialgeometrie im Gro\ss en}, 28 June - 4 July 2015, and contains no new results. 
\end{abstract}



\section{Introduction}	
 
 Higgs bundles were first studied by Nigel Hitchin in 1987, and appeared as solutions of Yang-Mills self-duality equations on a Riemann surface \cite{LPS_N1}. 
Classically, a {\it Higgs bundle} on a compact Riemann surface $\Sigma$ of genus $g\geq 2$ is a pair $(E,\Phi)$ where 
 $E$ is a holomorphic vector bundle  on $\Sigma$, 
and  $\Phi$, the {\it  Higgs field}, is a holomorphic 1-form in $H^{0}(\Sigma, {\rm End}_{0}(E)\otimes K_{\Sigma}),$
 for $K_\Sigma$ the cotangent bundle of $\Sigma$ and ${\rm End}_{0}(E)$ the traceless endomorphisms of $E$.  
Higgs bundles can also be defined for complex groups $G_{c}$, and through stability conditions, one can construct their moduli spaces $\mathcal{M}_{G_c}$. 

A natural way of studying the moduli space of Higgs bundles   is through the {\it Hitchin fibration}, 
  sending the class of a Higgs bundle $(E,\Phi)$ to the coefficients of the characteristic polynomial $\det(x I -\Phi )$. The generic fibre  is an abelian variety, which can be seen through line bundles on an algebraic curve $S$, the {\it spectral curve} associated to the Higgs field. The \textit{spectral data}  is then given by a line bundle on $S$ satisfying certain conditions, and it provides a geometric description of the fibres of the Hitchin fibration. For instance in the case of classical Higgs bundles, the smooth fibres can be seen through spectral data as Jacobian varieties of $S$.

   \section{Higgs bundles and branes}
  
 The smooth locus of the moduli space $\mathcal{M}_{G_c}$ of $G_c$-Higgs bundles on a compact Riemann surface $\Sigma$ for a  a complex reductive Lie group $G_c$ is a hyper-K\"ahler manifold, so there are natural complex structures $I,J,K$ obeying the same relations as the imaginary quaternions (adopting the notation of \cite{LPS_Kap}).
Adopting the language of physicists,  a Lagrangian submanifold of a symplectic manifold is called (the base of) an {\em A-brane} and a complex submanifold (the base of) a {\em B-brane}. A submanifold of a hyper-K\"ahler manifold may be of type $A$ or $B$ with respect to each of the complex or symplectic structures, and thus choosing a triple of structures one may speak of branes of type $(B,B,B), (B,A,A), (A,B,A)$ and $(A,A,B)$.  
 
  It is hence natural to construct  different families of branes inside the moduli space $\mathcal{M}_{G_c}$, as was first done in  \cite{LPS_Kap} (see also \cite{LPS_Witten}).  
    Together with D. Baraglia we introduced a naturally defined triple of commuting real structures $i_1,i_2,i_3$ on $\mathcal{M}_{G_c}$, and through the  {\it spectral data } gave a detailed picture of their fixed point sets as different types of branes. 
Given  $(E,\Phi)$   a $G_c$-Higgs bundle on $\Sigma$,   consider pairs $(\overline{\partial}_A , \Phi)$, where $\overline{\partial}_A$ denotes a $\overline{\partial}$-connection on $E$ defining a holomorphic structure, and $\Phi$ is a section of $\Omega^{1,0}(\Sigma , {\rm ad}(E))$, for ${\rm ad}(E)$ the adjoint bundle of $E$. Through the Cartan involution $\theta$ of a real form $G$ of $G_c$ one obtains 
\begin{equation}
i_1(\bar \partial_A, \Phi)=(\theta(\bar \partial_A),-\theta( \Phi)).
\end{equation}
  Moreover, a real structure  $f : \Sigma \to \Sigma$ on $\Sigma$  induces an involution $i_2$   given by 
\begin{equation}
i_2(\bar \partial_A, \Phi)=(f^*(  \partial_A),f^*( \Phi^*  ))= (f^*(\rho(\bar \partial_A)), -f^*( \rho(\Phi) )).
\end{equation} 
  Lastly, by looking at $i_3 = i_1 \circ i_2$, one obtains   
$i_3(\bar \partial_A, \Phi)=(f^* \sigma(\bar \partial_A),f^*\sigma( \Phi)).$
The fixed point sets of $i_1,i_2,i_3$ are branes of type $(B,A,A),(A,B,A)$ and $(A,A,B)$ respectively. In \cite{LPS_real}, with D. Baraglia we studied  these branes through the associated {\it spectral data} and described the topological invariants involved using $KO$, $KR$ and equivariant $K$-theory. In particular, it was shown that amongst the fixed points of $i_1$  are solutions to the Hitchin equations with holonomy in  $G$.
 
 \section{Higgs bundles and $(A,B,A)$-branes}
  Amongst the fixed points of the involution $i_2$ are representations of   $\pi_1(\Sigma)$ that extend to certain 3-manifolds $M$ whose boundary is $\Sigma$.    Indeed, consider the space  $\overline{\Sigma} = \Sigma \times [-1,1]$ with involution $\tau(x,t) = (f(x),-t)$, for $f$ the anti-holomorphic involution on $\Sigma$ giving $i_2$. The quotient $M = \overline{\Sigma}/\tau$ is a $3$-manifold with boundary $\partial M = \Sigma$, and satisfies the following:  
 
\begin{theorem}[Baraglia-Schaposnik~\cite{LPS_aba}]  Let $(E, \Phi)$ be a fixed point of $i_2$ with simple holonomy. Then the associated connection extends over $M$ as a flat connection if and only if the class $[E] \in \tilde{K}^0_{\mathbb{Z}_2}(\Sigma)$ in reduced equivariant $K$-theory is trivial.\end{theorem}
 
Since Langlands duality can be seen in terms of Higgs bundles as a duality between the fibres of the Hitchin fibrations for $\mathcal{M}_{G_c}$ and $\mathcal{M}_{^LG_{c}}$, for $^LG_c$ the Langlands dual group of $G_c$ (as was first seen in \cite{LPS_Tamas}), it is natural to ask what the duality between branes should be. In \cite{LPS_real} we proposed the following:

\begin{conjecture}[Baraglia-Schaposnik \cite{LPS_real}]
For $^L\rho$ the compact structure of $^LG_c$, the support of the dual brane of the fixed point set of $i_2$ is   the fixed point set in $\mathcal{M}_{^LG_c}$ of 
\begin{equation}
^Li_2(\bar \partial_A, \Phi)=  (f^*(^L\rho(\bar \partial_A)), -f^*( ^L\rho(\Phi) )).
\end{equation} 
\end{conjecture} 
 
 \section{ Higgs bundles and  $(B,A,A)$-branes }
 
Higgs bundles can be defined for complex Lie groups $G_c$, as well as for real forms $G$ of $G_c$. Moreover, as seen in \cite{LPS_Kap}, the moduli spaces of real Higgs bundles
$\mathcal{M}_{G}$ lie as $(B,A,A)$ branes inside the moduli spaces $\mathcal{M}_{G_c}$.  
 It is thus natural to ask how $\mathcal{M}_{G}$  intersects the Hitchin fibration for the complex moduli space $\mathcal{M}_{G_c}$, which is the main subject of \cite{LPS_yoU, LPS_thesis} (see also \cite{LPS_spectral}). Moreover, considering Langlands duality,  in \cite{LPS_real} we propose the following:
 
 \begin{conjecture}[Baraglia-Schaposnik \cite{LPS_real}]
 The support of the dual brane to the fixed point set of $i_1$ is the moduli space $\mathcal{M}_{\check{H}}\subset \mathcal{M}_{^LG_{c}}$ of $\check{H}$-Higgs bundles for $\check{H}$ the so-called Nadler group, the group associated to the Lie algebra $\check{\mathfrak{h}}$ in \cite[Table 1]{LPS_Nadler}. 
 \end{conjecture}

\subsection{On $(B,A,A)$-branes having finite intersection with smooth fibres}

 In the case of Higgs bundles for a split real form $G$ of a complex reductive Lie group $G_c$, from \cite{LPS_thesis} one has the following description of the intersection:
 
\begin{theorem}[Schaposnik \cite{LPS_thesis}]The moduli space $\mathcal{M}_{G}$ as sitting inside $\mathcal{M}_{G_c}$ is given by points of order two in the smooth fibres of the Hitchin fibration $h:\mathcal{M}_{G_c}\rightarrow \mathcal{A}_{G_c}$.\end{theorem}

This result is used in \cite{LPS_monodromy} to study the moduli space of $SL(2,\mathbb{R})$-Higgs bundles from a combinatorial point of view. The above theorem was also be used when studying $L$-twisted Higgs bundles $(E,\Phi)$ in \cite{LPS_mono2}, where the Higgs field is now twisted by any line bundle $L$ obtaining $\Phi:E\rightarrow E\otimes L$.

 \begin{theorem}[Baraglia-Schaposnik \cite{LPS_mono2}] The monodromy action of an element $\tilde s_\gamma$ in the Hitchin base is the automorphism of $ H^1(S , \mathbb{Z})$ induced by a Dehn twist of the spectral curve $S$ around a loop $l_\gamma$. Let $c_\gamma \in H^1(S , \mathbb{Z})$ be the Poincar\'e dual of the homology class of $l_\gamma$. Then the monodromy  of  $\tilde s_\gamma$ acts on $H^1( S , \mathbb{Z})$ as a Picard-Lefschetz transformation:
$x \mapsto x + \langle c_\gamma , x \rangle c_\gamma.$
 \end{theorem}
 As an application of this theorem, one can obtain information about some real character varieties $Rep(G)$. For instance, with D.~Baraglia we prove in \cite{LPS_mono2} that there are $3.2^{2g}+g-3$ connected components for $Rep(GL(2,\mathbb{R}))$,  and recover  the number of components of maximal representations $Rep_{2g-2}(Sp(4,\mathbb{R}))$,  given by $3.2^{2g}+2g-4$.
 
 \subsection{On $(B,A,A)$-branes having no intersection with smooth fibres}

For other real forms the brane $\mathcal{M}_{G}$ may lie   inside the singular fibres of the Hitchin fibration for $G_c$-Higgs bundles. By extending the approach from \cite{LPS_thesis} for $Sp(2p,2p)$-Higgs bundles,   together with N. Hitchin we showed in \cite{LPS_nonabelian} that this situation  appears naturally when considering Higgs bundles   corresponding  to flat connections on $\Sigma$ with holonomy in the real Lie groups $G=SL(m,\mathbb{H})$ and $SO(2n,\mathbb{H})$ (i.e.,  $SU^*(2m)$ and $SO^*(2n)$). 
  \vspace{0.1 in }

 \begin{theorem}[Hitchin-Schaposnik \cite{LPS_nonabelian}] The fibres of the $(B,A,A)$-brane in the Hitchin fibration for $SL(m,\mathbb{H})$, $SO(2n,\mathbb{H})$ and $Sp(2m,2m)$-Higgs bundles  are not abelian varieties, but are instead moduli spaces of rank 2 bundles on a spectral curve, satisfying certain natural stability conditions.\end{theorem}

\section*{Acknowledgements}  
The author would like to thank the organisers of  Oberwolfach's workshop {\it Differentialgeometrie im Gro\ss en}, 28 June - 4 July 2015,  for providing an ideal environment for collaboration.
The work of LPS is partially supported by NSF grant DMS-1509693.

\end{document}